\documentclass[12pt]{article}
\usepackage{amsmath,amsthm,amsfonts,amssymb}
\newtheorem{thm}{Theorem}[section]
\newtheorem{lemma}[thm]{Lemma}
\newtheorem{defn}[thm]{Definition}
\newtheorem{cor}[thm]{Corollary}

\newtheorem{prop}[thm]{Proposition}

\theoremstyle{definition}

\newtheorem{com}[thm]{Remark}
\newtheorem{ques}[thm]{Question}
\newtheorem{ex}[thm]{Example}

\newtheorem{remark}[thm]{Remark}
\theoremstyle{remark}
\DeclareMathOperator{\Fix}{Fix}
\DeclareMathOperator{\Diff}{Diff}
\DeclareMathOperator{\MCG}{MCG}
\DeclareMathOperator{\Ham}{Ham}
\DeclareMathOperator{\Homeo}{Homeo}

\DeclareMathOperator{\Per}{Per}
\DeclareMathOperator{\supp}{supp}

\DeclareMathOperator{\egr}{egr}
\DeclareMathOperator{\ax}{Ax}
\newcommand{\R}{\mathbb R}
\newcommand{\RP}{{\mathbb R \mathbb P}}

\newcommand{\N}{{\cal N}}

\newcommand{\T}{\mathbb T}
\newcommand{\Z}{\mathbb Z}
\title{Distortion Elements in Group actions on surfaces}
\author{John Franks\thanks{Supported in part by NSF grant DMS0099640.}\ \ 
and Michael Handel\thanks{Supported in part by NSF grant DMS0103435.}}

\def\ti{\tilde}
\def\sinfty{S_{\infty}}
\def\sl3z{SL(3, \mathbb Z)}

\def\D{{\cal D}}
\def\G{{\cal G}}
\def\H{{\cal H}}

\def\Q{{\mathbb Q}}

\begin{document}
\maketitle
\begin{abstract}
If $\G$ is a finitely generated group with generators
$\{g_1,\dots,g_j\}$ then an infinite order element $f \in \G$ is a
{\em distortion element} of $\G$ provided $\displaystyle{\liminf_{n
\to \infty} |f^n|/n = 0,}$ where $|f^n|$ is the word length of $f^n$
in the generators.  Let $S$ be a closed orientable surface and let
$\Diff(S)_0$ denote the identity component of the group of $C^1$
diffeomorphisms of $S$.  Our main result shows that if $S$ has genus
at least two and if $f$ is a distortion element in some finitely
generated subgroup of $\Diff(S)_0$, then $\supp(\mu) \subset \Fix(f)$
for every $f$-invariant Borel probability measure $\mu$. Related
results are proved for $S = T^2$ or $S^2$.

For $\mu$ a Borel probability measure on $S$, denote the group of
$C^1$ diffeomorphisms that preserve $\mu$ by $\Diff_{\mu}(S)$.  We
give several applications of our main result showing that certain
groups, including a large class of higher rank lattices, admit no
homomorphisms to $\Diff_{\mu}(S)$ with infinite image.
\item

\end{abstract}

\section{Introduction and Notation}

In this paper $S$ is a closed orientable surface and $\mu$ is a, not
necessarily smooth, Borel probability measure on $S$. We denote the
group of $C^1$ diffeomorphisms [preserving $\mu$] by $\Diff(S)$ [$\Diff_{\mu}(S)$] and its
identity component by $\Diff(S)_0$  [$\Diff_{\mu}(S)_0$].  The support of $\mu$ is
denoted $\supp(\mu)$.  The set of fixed points for $f$ is denoted $\Fix(f)$ and the set of periodic points for $f$ is denoted $\Per(f)$.  

\begin{defn}
If $\G$ is a finitely generated group with generators
$\{g_1,\dots,g_j\}$ then $f \in \G$ is said to be a {\em distortion element} of
$\G$ provided $f$ has infinite order and
\[
\liminf_{n \to \infty} \frac{|f^n|}{n} = 0,
\]
where $|f^n|$ is the word length of $f^n$ in the generators
$\{g_1,\dots,g_j\}$.  
If $\G$ is not finitely generated, then we say that $f \in \G$ is distorted in $\G$ if it is distorted in some finitely generated subgroup of $\G$.
\end{defn}

Equivalently, $f$ is distorted if the {\em translation length} of $f$
in $\G$ as defined in \cite{gs} is zero.

\begin{remark}\label{rem: distortion}
It is straightforward to see that the property of being a distortion
element is independent of the set of generators used to measure
$|f^n|.$  It is also immediate that if $f$ is a distortion element
then so is $f^k$ for all $k > 0.$  The structure theorem for finitely
generated abelian groups implies that they have no distortion elements and
hence no abelian group contains any distortion elements.
\end{remark}

It is easy to check (see Lemma~\ref{lem: dist_heis}) that central
elements of the three dimensional Heisenberg group are distortion
elements.  Lubotzky, Mozes and Raghunathan \cite{lmr} proved that
irreducible non-uniform lattices in higher rank Lie groups have distortion
elements.  (They prove the stronger result that there are elements in
the group whose word length has logarithmic growth; in this paper we
use only that the growth is sublinear.) On the other hand, there are
no distortion elements in biautomatic groups \cite{gs}, mapping class
groups \cite{flm} and the outer automorphism group of a finitely
generated free group \cite{a}.  In particular, none of these groups
have a subgroup that is isomorphic to the three dimensional Heisenberg
group.

Polterovich (1.6.C of \cite{P}) proved that if $S$ is a closed surface of genus at least two and $\omega$ is a smooth measure on $S$, then  finitely generated subgroups of $\Ham_{\omega}(S)$, the group of Hamiltonian diffeomorphisms on $S$,  do not contain distortion elements.  There is also a closely related result (Theorem 1.6.B of \cite{P}) for finitely generated subgroups of $\Diff_{\omega}(S)$.

Our main theorem is a generalization of these results from \cite{P}.
We do not assume that the invariant measure is smooth and in
particular do not assume that its support is all of $S$.
In several respects there are parallels between properties of 
$\Diff(S^1)$ and $\Diff_\mu(S)$ for $S$ a surface. One of the striking
properties of homeomorphisms of $S^1$ is that
periodic points all have the same period and 
if there are periodic points they contain the support of any invariant
measure.  Of course, in general this is far from true for surfaces.
However, our main result shows this property does generalize for distortion
elements in $\Diff_\mu(S)$.

\begin{thm}  \label{mainDistort} 
Suppose that $S$ is a closed oriented surface, that $f$ is a distortion element in $\Diff(S)_0$ and that $\mu$ is an $f$-invariant Borel probability measure.
\begin{enumerate}
\item
If $S$ has genus at least two then $\supp(\mu) \subset \Fix(f)$. 
\item
If $S = T^2$ and $\Per(f) \ne \emptyset$ then $\supp(\mu)
\subset \Fix(f^n),$  where $n > 0$ is the
minimal value for which $\Fix(f^n) \ne \emptyset$.
\item
If $S = S^2$ and if $f^n$ has at least three fixed points
for some smallest $n>0$, then 
$\supp(\mu) \subset \Fix(f^n).$
\end{enumerate}
\end{thm}

\begin{remark}   Every periodic orbit is the support of some invariant measure.  Theorem~\ref{mainDistort} therefore implies that if $f$ is a distortion element in $\Diff(S)_0$   and if $S$ has genus at least two then $\Per(f) = \Fix(f)$.  Similarly, if $S$ is $T^2$ or $S^2$ and if $n> 0$ is as in (2) or (3) above, then  $\Per(f) = \Fix(f^n)$.  
\end{remark}

We do not know if the hypothesis that $\Per(f)$ contains at least three points is necessary in (3).

This result is motivated by an interest in understanding the 
algebraic properties of $\Diff_{\mu}(S)_0$.  We will apply
it to obtain results showing that certain finitely generated
groups can only have homomorphic images in $\Diff_{\mu}(S)_0$
which are finite.

Recall that a group $\G$ is called {\em almost simple} provided
every normal subgroup is finite or has finite index.  Also $\G $
is called {\em Kazhdan} or said to satisfy {\em property T} provided
the identity is isolated in the unitary dual $\hat \G$ (see
\cite{K} or \cite{Z1}).

\begin{thm}  \label{thm: finite_image} 
Suppose that $S$ is a closed oriented surface of positive genus
equipped with a Borel probability measure $\mu$ 
and that $\G$ is a finitely generated group which is almost simple and
possesses a distortion element $u$.  Suppose further that either $\mu$ has infinite support or that $\G$ is a Kazhdan group.  Then any homomorphism $\phi: \G
\to \Diff_{\mu}(S)$ has finite image.  The result is valid for $S =
S^2$ with the additional assumption that $\phi(u)$ has at least three
periodic points.\end{thm}

In the case $S = S^2$ we may remove the rather artificial hypothesis on
periodic points of $\phi(u)$ at the expense of
replacing the hypothesis of the
existence of a distortion element with the stronger hypothesis that
$\G$ contains a subgroup isomorphic to the three-dimensional
Heisenberg group.  In this case the additional assumptions on $\mu$ are not necessary.

\begin{thm}  \label{thm: Heisenberg} 
Suppose that $S$ is a closed oriented surface
equipped with a Borel probability measure $\mu$, and that $\G$ is a
finitely generated almost simple group which has a subgroup $\H$
isomorphic to the three-dimensional Heisenberg group.
Then any homomorphism $\phi: \G \to \Diff_{\mu}(S)$ has finite image.
\end{thm}

There are two almost immediate corollaries obtained by applying these
results to higher rank lattices.  Their proofs amount to checking that
the hypotheses on the lattice are sufficient to give the hypotheses of
one of the results above.  We defer these proofs to the last section.

\begin{cor} \label{cor: lattice}  
Let $\G$ be an irreducible non-uniform lattice in a semisimple real
Lie group of real rank at least two.  Assume that the Lie group is
connected, without compact factors and with finite center.  Suppose that  $S$ is
a closed oriented surface of genus at least one equipped with   a Borel
probability measure $\mu$. Then every homomorphism
$\phi : \G \to \Diff_{\mu}(S)$ has finite image.
\end{cor}

In case that $\mu$ is smooth, 
Corollary~\ref{cor: lattice} is due to Polterovich \cite{P}.  There are
also results due to Farb and Shalen \cite{FS} and Ghys \cite{G}
which show that large classes of lattices admit no analytic actions
on surfaces except those which factor through a finite group.  No
invariant measure is required, but the action must be analytic.

Any finite index subgroup of $SL(n,\Z), \ n >2,$ contains a subgroup isomorphic
to the Heisenberg group and  is also almost simple, 
(as a result of the Margulis normal subgroups theorem; see \cite{Mar}).
Therefore as an immediate consequence of Theorem~(\ref{thm: Heisenberg})
we have the following.

\begin{cor} \label{cor: Heisenberg}  
Suppose that $\Gamma$ is a finite index subgroup of $SL(n,\Z) $.  If
$\mu$ is a Borel probability measure on a closed oriented
surface $S$ and $n \ge 3$,  
then every homomorphism $\phi : \Gamma \to \Diff_{\mu}(S)$ has finite 
image.
\end{cor}

It is known that some interesting classes of groups have the property
that their commutator subgroup consists entirely of distortion
elements (see \cite{BM}, Corollary (1.4) for example).
With this hypothesis our main result gives another interesting corollary.

\begin{cor} \label{cor: commutator}  
Suppose that $S$ is a closed oriented surface of genus at least one
equipped with a Borel probability measure $\mu$ whose support is all
of $S$. Suppose also that $\Gamma$ is a finitely generated group such
that $[\Gamma,\Gamma]$ is generated by distortion elements.  Then
every homomorphism $\phi : \Gamma \to \Diff_{\mu}(S)$ has an
abelian image.
\end{cor}

In the same vein we have the following.

\begin{cor} \label{cor: nilpotent}  
Suppose that $S$ is a closed oriented surface equipped with a Borel
probability measure $\mu$ whose support is all of $S$ and that $\N \subset \Diff_{\mu}(S)_0$ is nilpotent.  If $S \ne S^2$ then   $\N $ is abelian.  If $S = S^2$ then $\N$ has an index two abelian subgroup.
\end{cor}

To conclude this section we give several illustrative examples.
An example of G. Mess shows that at least some translations on
$\T^2$ are distorted.

\begin{ex}[G. Mess]\label{Mess}
Suppose $\mu$ is Lebesgue measure on $\T^2.$ 
In the subgroup of $\Diff_\mu(\T^2)$ generated by the automorphism
\[
A = \begin{pmatrix}
2 & 1\\
1 & 1\\
\end{pmatrix}
\]
and a translation $T(x) = x + w$ where $w \ne 0$ is parallel to the 
unstable manifold of $A$, the element $T$ is distorted.
\end{ex}

\begin{proof}
Let $\lambda$ be the expanding eigenvalue of $A$.
The element 
$h_n = A^n T A^{-n}$ satisfies $h_n(x) = x + \lambda^n w$
and $g_n = A^{-n} T A^n$ satisfies $g_n(x) = x + \lambda^{-n} w$.
Hence $g_n h_n(x) = x +  (\lambda^n + \lambda^{-n}) w.$
Since $tr A^n = \lambda^n + \lambda^{-n}$ is an integer we conclude
\[
T^{tr A^n} = g_n h_n, \text{ so } |T^{tr A^n}| \le  4n +2.
\]
But
\[
\lim_{n \to \infty} \frac{n}{tr A^n} = 0,
\]
so $T$ is distorted.
\end{proof}

A similar example illustrating Theorem~\ref{mainDistort} starts
with an action on $S^1.$

\begin{ex}
Let $\G$ be the subgroup of $PSL(2,\Z[\sqrt{2}])$ generated by 
\[
A =
\begin{pmatrix}
\sqrt{2} -1 &  0\\
0 &  \sqrt{2}+1\\
\end{pmatrix}
\text{ and }
B =
\begin{pmatrix}
1 &  1\\
0 & 1 \\
\end{pmatrix}.
\]

Let $\lambda = \sqrt{2} +1 $ and note $\lambda^{-1} = \sqrt{2} -1$.
These matrices satisfy
\[
A^{-n}BA^n = 
\begin{pmatrix}
1 & \lambda^{2n} \\
0 & 1 \\
\end{pmatrix}
\]
and
\[
A^n BA^{-n} = 
\begin{pmatrix}
1 & \lambda^{-2n} \\
0 & 1 \\
\end{pmatrix}.
\]
It is easy to see that $m =\lambda^{2n} + \lambda^{-2n}$ is
an integer. Hence 
\[
(A^{-n}BA^n) (A^n BA^{-n}) =
\begin{pmatrix}
1 & \lambda^{2n} + \lambda^{-2n} \\
0 & 1 \\
\end{pmatrix}
= B^m.
\]
Thus $m$ grows exponentially in $n$ while $|B^m| \le 4n+2$ so $B$ is
distorted.  The group $\G$ acts naturally on $\RP^1$ (the lines in
$\R^2$ through the origin) which is diffeomorphic to $S^1$.  The
element $B$ has a single fixed point, the $x-$axis, and the only $B$
invariant measure is supported on this point.

In example 1.6.K of \cite{P} Polterovich considers the
embedding $\psi: \G \to PSL(2,\R) \times PSL(2,\R)$ where
$\psi(g) = (g, \bar g)$ with $\bar g$ being the conjugate of $g$
obtained by replacing an entry $a+b\sqrt{2}$ with $a-b\sqrt{2}.$
He points out that the image of $\psi$ is an irreducible
non-uniform lattice in a Lie group of real rank $2.$  Of course
$(B, \bar B) = (B, B)$ is a distortion element in $\psi(\G)$ and in 
the product action of $PSL(2,\R) \times PSL(2,\R)$ on $T^2 = S^1 \times S^1$ 
it has only one fixed point $(p, p)$ where $p$ is the fixed point of 
$B$ acting on $S^1.$  It is also clear that the only $(B, \bar B)$ invariant
measure is supported on this point.  It is easy to see that there are
elements of $\psi(\G)$ which do not fix this point, and hence there is
no probability measure invariant under all of $\psi(\G).$ 

\end{ex}

D. Calegari showed us an example which illustrates that 
$\Diff_\mu(S)$ 
cannot be replaced with $\Homeo_\mu(S)$  in Theorem~\ref{thm: Heisenberg}.

\begin{ex}[D. Calegari]\label{Calegari}
There is a $C^0$ action of the Heisenberg group on $S^2$ 
whose center is generated by an irrational
rotation.  Hence an irrational rotation of $S^2$ is distorted
in $\Homeo_\mu(S^2)$ where $\mu$ is a measure supported on the fixed
points of this rotation.
\end{ex}

\begin{proof}
Consider the linear homeomorphism of $\R^2$ given by 
$G(x,y) = (x+y, y)$ 
and the translation $H(x,y) = (x ,y +1)$.  We compute
$F = [G,H]$ to be a translation $F(x,y) = (x+1, y).$
This defines an action of the Heisenberg group on $\R^2$.
Let $C$ be the cylinder obtained by quotienting by the 
relation $(x ,y ) \sim (x + \alpha, y)$ for some 
$\alpha \in \R \setminus \Q$. The quotient
action is well defined.  The two ends of $C$ are fixed by every element
of the action and hence if we compactify $C$ to obtain
$S^2$ by adding a point at each end, we obtain an action
of the Heisenberg group by homeomorphisms on $S^2.$ 
\end{proof}

\begin{ques}
Is an irrational rotation of $S^1$ distorted in $\Diff(S^1)$ or
$\Homeo(S^1)?$  
Is an irrational translation distorted in $\Diff(T^2)_0?$ 
Example~\ref{Mess} shows that at least some irrational translations
are distorted $\Diff(T^2)$, but the construction requires generators
which are not isotopic to the identity.
Is an irrational rotation of $S^2$ distorted in
$\Diff(S^2)$ or area preserving diffeomorphisms of $S^2?$ 
Example~\ref{Calegari} shows that any irrational rotation of
$S^2$ is distorted in $\Homeo(S^2).$  
\end{ques}

Conversations with K. Parwani were valuable in the preparation of this
article.

\section{Rotation vectors} \label{rotation vectors}

Suppose that $S \ne S^2$, that $g \in \Diff(S)_0$ and that $\ti g :
\ti S \to \ti S$ is a lift to the universal cover that commutes with
all covering translations of $\ti S$.  If $S \ne T^2$, then there is a
unique such lift but if $S= T^2$ then every lift has this property.

Given $x \in S$, choose a lift $\ti x \in \ti S$.  For each $n > 0$,
choose a lift $\ti x_n$ of $x$ whose distance from $\ti g^n(\ti x)$ is
uniformly bounded, say by the diameter of $S$.  There is a covering
translation $T_n$ such that $T_n(\ti x) = \ti x_n$ and we denote the
element of $H_1(S)$ corresponding to $T_n$ by $\delta_n(x)$.  We
define the {\em rotation vector} of $x \in S$ with respect to $\ti g$
by
\[
\rho(x, \ti g) 
= \lim_{n \to \infty} \frac{1}{n}  \delta_n(x)) 
\in H_1(S, \R)
\]
if this limit exists.  It is clear that $\rho(x, \ti g)$ is
independent of the exact choice of $T_n$.  Moreover, if $A$ is a
covering translation and if $\ti x$ is replaced by $A(\ti x)$, then
$T_n$ can be replaced by $A T_n A^{-1}$. Thus $\rho(x, \ti g)$ depends only
on $x$ and $\ti g$ and not on the lift $\ti x$.

The following results are proved in \cite{Fr5}. 

If $g \in \Diff_{\mu}(S)_0$ then $\rho(x, \ti g)$ is defined for
almost every point with respect to $\mu$.  If $S \ne T^2$
that the {\it mean rotation vector} $\rho_{\mu}(\ti g) := \int \rho(x,
\ti g)d\mu : \Diff_{\mu}(S)_0 \to H_1(S, \R)$ is a homomorphism.

There is a natural identification of $H_1(T^2,\R)$ with the universal
cover $\R^2$ of $T^2$.  Moreover, if $T : \R^2 \to \R^2$ is a covering
translation then $\rho(x, T\ti f) = T \rho(x, \ti f)$.  It follows
that the image $\hat \rho(x,f) \in T^2$ of $\rho(x,\ti f)$ is
independent of the choice of lift $\ti f$.  There is an induced
homomorphism $\hat \rho_{\mu} : \Diff_{\mu}(T)_0 \to T^2$.  Abusing
terminology slightly we will also refer to $\hat \rho_{\mu}(f)$
as a rotation vector.

These definitions works equally well for a homeomorphism $g : A \to A$
of a closed annulus.  In that case, $H_1(A,\R)$ is identified with
$\R$, $\rho(x,\ti g)$ is called the {\it rotation number} of $x$ with
respect to $\ti g$ and $\rho_{\mu}(\ti g) $ is called the {\it mean
rotation number} of $\mu$ with respect to $\ti g$.  There is an induced
homomorphism $\hat \rho_{\mu}: \Homeo_{\mu}(A)_0 \to S^1$.

\begin{com}  \label{zero mean rotation}
Suppose that $g \in \Diff_{\mu}(S)_0$ or $g \in \Diff_{\mu}(A)_0$ .
Then $\hat \rho_{\mu}( g) = 0$ if and only if there is a lift $\ti g$
such that $\rho_{\mu}(\ti g)=0$.
\end{com} 

\begin{lemma}  \label{ergodic measure}  Suppose that 
$M = T^2$ or $M = A$, that $g \in \Homeo(M)_0$ and that 
for some lift $\ti g$to the universal cover $\rho(x, \ti g) = 0$ 
for all $x$ in the interior of $M$ for which it is defined.
Then
\begin{enumerate}
\item $\Fix(g) \ne \emptyset$.
\item If in the case $M=A$ there is an invariant measure $\mu$ 
with $\mu( A \setminus \partial A) >0$, then $g$ has a fixed point in the
interior of $A$.
\end{enumerate}
\end{lemma}

\proof 
Suppose first that $M = T^2.$
It is an immediate
consequence of the hypotheses and the definitions that 
$\rho_{\mu}(\ti g) = 0$ for every invariant measure $\mu$.
There is always an invariant ergodic measure $\mu$.  The lemma
therefore follows from applying Theorem 3.5 of \cite{Fr}.

For the case $M = A$ we note that under the hypotheses of (2) there
is an invariant ergodic measure $\nu$  with $\nu( A \setminus \partial A) >0$.
It then must be the case that $\nu(\partial A) = 0$ and 
that $\rho_{\nu}(\ti g) = 0$.  The result then follows
from the following proposition.  
\endproof

\begin{prop}  \label{prop: ergodic_annulus}  
Suppose that $g: A \to A$ is a homeomorphism of the annulus which
is isotopic to the identity and possesses an ergodic invariant measure
$\nu$ whose support intersects the interior of $A$.  
If $\rho_{\nu}(\ti g) = 0$ for some lift $\ti g: \ti A \to \ti A$ of
$g$ to the universal cover $\ti A$ then $\ti g$ has a fixed point
in the interior of $\ti A.$
\end{prop}

\begin{proof}  
Let $P$ be a disk in the interior of $A$ with $\nu(P) >0$ and let $h$
be the first return map of $g$ on $P$.  

Denote the full pre-image of
$P$ in $\ti A$ by $\ti P$ and the lift of $h$ corresponding to 
$\ti g$ by $\ti h$.  Projection
$p_1 :A \to S^1$ onto the first coordinate lifts to a projection $\ti
p_1: \ti A \to \R$.  For $x$ in the domain of $h$, define $\psi_0(x) =
(\ti p_1(\ti g(\ti x)) - \ti p_1(\ti x))$ where $\ti x$ is a lift of $x$
and define $\psi(x) =(\ti p_1(\ti h(\ti x)) - \ti p_1(\ti x))$.
From the fact that  
\[
\int_A \psi_0 \ d\nu = \rho_{\nu}(\ti g) = 0
\]
it is straightforward to check that
\[
\int_P \psi\ d\nu = 0.
\]

Using this, Proposition (12.1) from \cite{fh:periodic} asserts that if
\[
S(n,x) = \sum_{i=0}^{n-1} \psi( h^i(x))
\]
then for a full measure subset of $P$ there are infinitely many
values of $n >0$ with $|S(n,x)| < 1.$
Hence $p_1(h^n(\ti x)) - p_1(\ti x) <1$ for infinitely many values of $n$.
We conclude that $\ti g$ has an interior recurrent point and by the
Brouwer plane translation theorem has an interior fixed point.
\end{proof}

\section{Hyperbolic Geometry} \label{hyperbolic}

Some of our proofs rely on mapping class group techniques that use
hyperbolic geometry.  In this section we establish notation and recall
standard results about hyperbolic structures on surfaces. 

 Let $S$ be
a closed orientable surface.  We will say that a connected open subset
$M$ of $S$ has {\em negative Euler characteristic} if $H_1(M, \R)$ is
infinite dimensional or if $M$ is of finite type and the usual
definition of Euler characteristic has a negative value. All such $M$ have complete hyperbolic structures; further details can be found, for example, in  \cite{fh:periodic}.

We use the Poincar\'e disk model for the hyperbolic plane $H$.  In
this model, $H$ is identified with the interior of the unit disk and
geodesics are segments of Euclidean circles and straight lines that
meet the boundary in right angles. A choice of hyperbolic structure on
$M$ provides an identification of the universal cover $\ti M$ of $M$
with $H$.  Under this identification covering translations become
isometries of $H$ and geodesics in $M$ lift to geodesics in $H$.  The
compactification of the interior of the unit disk by the unit circle
induces a compactification of $H$ by the \lq circle at infinity\rq\
$\sinfty$.  Geodesics in $H$ have unique endpoints on $\sinfty$.

 Each covering translation $T$ extends to a homeomorphism (also called) $T : H \cup \sinfty \to H \cup \sinfty$. The fixed point set of a non-trivial $T$ is either one or two points in $\sinfty$. We denote these point(s) by $T^+$ and $T^-$, allowing the possibility that $T^+ = T^-$.  If $T^+ = T^-$, then $T$ is said to be {\it parabolic}.  If $T^+$ and $T^-$ are distinct, then    $T$ is said to be {\it hyperbolic} and we may assume that $T^+$ is a sink and $T^-$ is a source. In this case the geodesic connecting $T^-$ to $T^+$ is said to be the {\it axis of $T$}.

Suppose now that $f \in \Diff(S)$ and  that $K \subset \Fix(f)$. Each component  of $S \setminus \Fix(f)$ is $f$-invariant by \cite{brn-kis}.   If $M$ is a component of $S \setminus K$ with negative Euler characteristic, then we   denote $f|_M$ by $h$ to keep a clear distinction between $f$ and $f|_M$ and to simplify notation. 
We use the identification of $H$ with $\ti M$ and write $\ti h : H \to
H$ for lifts of $h: M \to M$ to the universal cover.  A fundamental
result of Nielsen theory is that every lift $\ti h  : H \to H$ extends
uniquely to a homeomorphism (also called) $\ti h   : H \cup \sinfty \to
H \cup \sinfty$.  (A proof of this fact appears in Proposition 3.1 of
\cite{han:fpt}).  If $h   : M \to M$ is isotopic to the identity then
there is a unique lift $\ti h  $, called the {\it identity lift,} that
commutes with all covering translations and whose extension over
$\sinfty$ is the identity.

In certain contexts it is convenient to blow up isolated ends of $M$. More precisely, we  compactify an isolated end   of $M$  by adding a circle via the following well known lemma,  a proof of which can be found in Lemma 5.1 of \cite{fh:periodic}.  This construction does not depend on the hyperbolic structure of $M$ and can be applied equally well when $M$ is topologically an annulus.

\begin{lemma}\label{extension}   
Suppose that $U \subset S$ is an open $f$-invariant set and that $E$ is an
isolated end of $U$ whose frontier is contained in Fix($f$).  Then $E$
can be compactified by a circle $C$  and $f|_U$ has a continuous extension $\overline {f|_C}$
  over $C$. If the frontier of $E$ is not a single point
then $ \overline {f|_C}=$ identity.  Moreover, if $\gamma$ is a smooth path in $U$ with an endpoint $p$ in the frontier of $U$, then $\gamma \setminus \{p\}$ extends to a path in $\bar U = U \cup C$ with an endpoint $P$ in $C$.
\end{lemma}

This \lq blowing up\rq\ is functorial in the sense that if $g: S \to
S$ preserves $U$ and pointwise fixes the frontier of $U$ then
$\overline {f|_C \circ g|_C} = \overline {f|_C} \circ \overline
{g|_C}.$ (This is immediate if the frontier of $E$ is not a single
point and follows easily from the smoothness of $f$ and $g$ in the
remaining case.) Similarly, if $\gamma_1$ and $\gamma_2$ are smooth
paths in $U$ with a common endpoint $p$ in the frontier of $U$, and
if $\gamma_2$ and $\gamma_2$
are homotopic relative to the frontier of $U$, then
the extensions of $\gamma_1 \setminus \{p\}$ and $\gamma_2 \setminus
\{p\}$ have a common endpoint $P$ in $C$.

If   finitely many   isolated ends of $M$  have been blown up,  then the resulting surface, which we continue to denote by $M$,  has a complete  hyperbolic structure  but now its universal
cover  is naturally identified with the intersection of $H$ with the  convex
hull of a Cantor set $C \subset \sinfty$.  The frontier of $\ti M$  in $H \cup \sinfty$ is
the union of $C$ with the full pre-image $\partial \ti M$ of $\partial M$.  If $h$ is
isotopic to the identity, then the identity lift commutes with all
covering translations and extends to a homeomorphism of $\partial \ti M \cup C$  that fixes $C$.

\section{Some results when $K = \Fix(f)$}

The main lemma of this section is a consequence of the results and methods of \cite{fh:periodic}.    We  begin with a pair of definitions.  The first is {\it linear tracing}.

\begin{defn} \label{tracing} 
Suppose that $f \in \Diff(S)$, that $M$ is a component of $S \setminus
\Fix(f)$ with negative Euler characteristic and that $\beta$ is an
essential simple closed curve in $M$ whose isotopy class is fixed by
$h=f|_M$. If $\beta$ is peripheral in $M$ assume that the end that it
encloses has been blown up to a boundary component.  Choose a covering
translation $T :H \cup \sinfty \to H \cup \sinfty$ whose axis $\ax_T$
projects to a simple closed curve that is isotopic to $\beta$.
Identify $\ax_T$ with $\R$ so that the action of $T$ on $\ax_T$
corresponds to unit translation on $\R$.  Let $\ti p: H \cup (\sinfty
\setminus\{T^{\pm}\}) \to \R$ be a $T$-equivariant projection of $ H
\cup \sinfty$ onto $\ax_T$ (e.g. orthogonal projection) followed by
the identification of $\ax_T$ with $\R$.  Let $\ti h :H \cup \sinfty
\to H \cup \sinfty $ be a lift of $h$ that commutes with $T$.  We say
that $x \in M$ {\it linearly traces $\beta$} if there is a lift $\ti x
\in H$ such that
\[
\liminf_{n \to \infty} \frac{| \ti p(\ti h^n(\ti x))- \ti p(\ti
 x)|}{n} > 0.
\]
\end{defn}

Since any two $T$-equivariant projections of $H \cup \sinfty$ onto $\ax_T$ differ by a bounded amount, Definition~\ref{tracing} is independent of the choice of projection.  We make use of this flexibility in Lemma~\ref{spread and tracing}. 

The surface $S$ is assumed to have an underlying Riemannian metric so
a smooth closed curve $\tau \subset S$ has a well defined length $l_S(\tau)$.
Define the {\it exponential growth rate} by
$$
\egr(f,\tau) = \liminf_{n \to \infty}\frac{\log(l_S(f^n(\tau))}{n}
$$.

We can now state and prove the main result of this section.

\begin{lemma} \label{fromPeriodic} Suppose that $f \in \Diff_{\mu}(S)_0$ has infinite order and that $\supp(\mu) \not \subset \Fix(f)$.  Suppose further that if $S = T^2$  then $\Fix(f) \ne \emptyset$ and that if $S= S^2$ then $\Fix(f)$ contains at least three points.  Then, after possibly replacing $f$ with an iterate,   at least one of the following holds:
\begin{description}
\item [(1)]There is a closed curve $\tau$ such that $\egr(f,\tau) > 0$.
\item [(2)] $f$ is isotopic relative to $\Fix(f)$ to a composition of non-trivial Dehn twists about a finite collection of non-peripheral, non-parallel, disjoint simple closed curves in $S \setminus \Fix(f)$.
 \end{description}
For the following, $f$ is isotopic to the identity relative to $\Fix(f)$.
 \begin{description}
\item [(3)]There is an $f$-invariant annular component $U$ of $S \setminus \Fix(f)$ with canonical extension $\overline{f|_U}: \bar U \to \bar U$ and there is $x \in U$ such that the rotation number of $x$ with respect to the lift of $\overline{f|_U}$ that fixes points on the boundary is non-zero.
\item [(4)] There exists a component $M$ of $S \setminus \Fix(f)$ with negative Euler characteristic, a simple closed curve $\beta \subset M$ that is essential in $M$ and there exists  $x \in M$ that linearly traces $\beta$ in $M$.  
\item [(5)] There exists $x \in S$ and a lift $\ti f$ such that $\Fix(\ti f) \ne \emptyset$ and such that the rotation vector $\rho(x,\ti f)$ is not zero. In this case $S \ne S^2$ and  $\ti f$ is the identity lift if $S \ne T^2$.  
\end{description}
\end{lemma}

\proof  We follow the proof of Theorem 1.1 of \cite{fh:periodic}.  

  Theorem 1.2 of \cite{fh:periodic} implies that if neither (1) nor (2) holds then $f$ is isotopic to the identity relative to $\Fix(f)$.   Applying this to iterates of $f$ we may   assume that $f^k$ is isotopic to the identity relative to $\Fix(f^k)$ for all $k \ge 1$.

  Each component $U$ of $S \setminus \Fix(f)$ is $f$-invariant by \cite{brn-kis}.  Choose a component $U$ of $S \setminus \Fix(f)$ with  $\mu(U) > 0$. Poincare recurrence implies that there are recurrent points in $U$.   The Brouwer plane translation theorem therefore implies that $U$ is not a disk.  If $U$ is an open annulus, then  by Lemma~\ref{extension}, $U$ can be compactified to a closed
annulus $\bar U$ in such a way that $f$ extends to the identity on at
least one boundary component.  Lemma~\ref{ergodic measure} implies that  (3) is satisfied.

We may therefore assume that $U$ has negative Euler characteristic; for notational consistency we denote $U$ by $M$. Lemma 3.7 of \cite{fh:group}, applied to the map obtained by redefining $f$ to be the identity on all components of $S \setminus \Fix(f)$ other than $M$, implies that $\Per(f|M) = \Fix(f|M)$.
 
The three preceding paragraphs replace the first four paragraphs of the proof of Theorem 1.1 of \cite{fh:periodic}.  The remainder of that proof, save the very last sentence, apply to our current context without change  and prove that either (4) or (5) hold.
\endproof

\section{Spread} \label{s: spread}

   Suppose  that $f \in \Diff(S)_0$, that $\gamma \subset S$ is a smoothly embedded path with distinct endpoints in $\Fix(f)$ and that $\beta$ is a simple closed curve that crosses $\gamma$ exactly once. We want to measure the relative rate at which points move \lq across
$\gamma$ in the direction of $\beta$\rq.  

 Let $A$ be the endpoint set of $\gamma$ and let $M$   be the surface with boundary obtained from $S$ by blowing up both points of $A$. We now think of $\gamma$ as a path in $M$ and of $\beta$ as a simple closed curve in $M$.  Assume at first that $S \ne S^2$ and  that  $M$ is equipped with a hyperbolic structure as in  section~\ref{hyperbolic}.  
 Choose non-disjoint extended lifts $\ti \beta \subset H \cup \sinfty$ and $\ti \gamma \subset H \cup \sinfty$  and let $T : H \cup \sinfty \to H \cup \sinfty$ be the covering translation corresponding to $\ti \beta$, i.e. $T^{\pm}$ are the endpoints of $\ti \beta$.    Denote $T^i(\ti \gamma)$ by $\ti \gamma_i$. Each $\ti \gamma_i$ is an embedded path in $H \cup \sinfty$ that intersects $\sinfty$ exactly in its endpoints.  Moreover, $\ti \gamma_i$ separates $\ti \gamma_{i-1}$ from $\ti \gamma_{i+1}$. 
 
An embedded smooth path $\alpha \subset S$ whose interior is disjoint from $A$ can be thought of as a path in $M$. For each lift $\ti \alpha \subset  H \cup \sinfty$, there exist $a < b$ such that $\ti \alpha \cap \ti \gamma_i \ne \emptyset$ if and only if $a < i < b$.  Define 
$$
\ti L_{\ti \beta, \ti \gamma}(\ti \alpha) = \max\{0,b-a-2\}
$$
and  
$$
L_{\beta,\gamma}(\alpha) = \max\{\ti L_{\ti \beta,\ti \gamma}(\ti \alpha)\}
$$
 as $\ti \alpha$ varies over all lifts of $\alpha$.  

Suppose now that $S= S^2$ and hence that $M$ is the closed annulus.  In this case $\ti M$ is identified with $\R \times [0,1]$, $T(x,y) = (x+1,y)$ and $\ti \gamma$ is an arc with endpoints in both components of $\ti M$.  With these modifications, $L_{\beta,\gamma}(\alpha)$ is defined as in the $S \ne S^2$ case.  

There is an equivalent definition of $L_{\beta,\gamma}(\alpha)$ that does not involve covers or blowing up.  Namely, $L_{\beta,\gamma}(\alpha)$  is the maximum value $k$ for which there exist subarcs $\gamma_0 \subset \gamma$ and $\alpha_0 \subset \alpha$ such that $\gamma_0\alpha_0$ is a closed path that is freely homotopic relative to $A$ to $\beta^k$. We allow the possibility that $\gamma$ and $\alpha$ share one or both endpoints.  The finiteness of $L_{\beta,\gamma}(\alpha)$ follows from the smoothness of the arcs $\alpha$ and $\gamma$.

\begin{com} \label{alg int} For future reference we observe that there is a subpath $\alpha_1 \subset \alpha_0$  so that the absolute value of the algebraic intersection of $\alpha_1$ with $\gamma$ is at least $L_{\beta,\gamma}(\alpha)-1$.  
\end{com}

Define the {\em spread} of
$\alpha$ with respect to $f, \beta$ and $\gamma$ to be 
\[
\sigma_{f, \beta,\gamma}(\alpha) 
=  \liminf_{n \to \infty} \frac{L_{\beta,\gamma}( f^n \circ \alpha)}{n}.
\]

\begin{com} \label{gamma independent} If $\gamma'$ is another smoothly embedded arc that crosses $\beta$ exactly once and that has the same endpoints as $\gamma$  then $\sigma_{f, \beta,\gamma}(\alpha)= \sigma_{f, \beta,\gamma'}(\alpha)$ for all $\alpha$.    This follows from the fact that  $\ti \gamma' $ is contained in the region bounded by $\ti \gamma_j$ and $\ti \gamma_{j+J}$ for some $j$ and  $J$ and hence $|L_{\beta,\gamma'}(\alpha) -L_{\beta,\gamma}(\alpha)| \le 2J$ for all $\alpha$.
\end{com}

\begin{lemma} \label{spread and tracing} Suppose that $M$ is a component of $S \setminus \Fix(f)$ with negative Euler characteristic, that $\beta$ is a simple closed essential curve in $M$ and that there exists a smoothly embedded  path $\gamma \subset S$ that has interior in $M$,  endpoint set $A \subset \Fix(f)$ and that  crosses $\beta$ exactly once.  If there exists $x \in M$ that linearly traces $\beta$ then  $\sigma_{f,\beta,\gamma}(\gamma) > 0$.
\end{lemma}

\proof Let $\ti \beta, \ti \gamma$ and $T$ be as in the definition of
spread.  Let $\ti p: H \cup (\sinfty \setminus\{T^{\pm}\}) \to \R$ and
$\ti x$ be as in Definition~\ref{tracing}.  We may assume that $\ti p
\ti \gamma = 0$ and $\ti p(\ti x) \in (-1,0)$.  Thus $\ti \gamma_k =
\ti p^{-1}(k)$ for all $k$.  For all sufficiently large $n$, $ \ti
p(\ti h^n(\ti x)) > Cn$ for some fixed constant $C > 0$.  
Since the endpoints of $\ti \gamma$ are fixed by $\ti h$ it follows
that there is a subpath $ \mu_n$ of $ \ti h^n(\ti \gamma)$ with
initial endpoint $z_1 \in \ti \gamma_1$ and terminal endpoint $z_2 \in
\ti \gamma_{Cn}$ and there is a subpath $\nu_n$ of $\ti \gamma_1$ with
terminal endpoint $z_1$ and initial endpoint $T^{1-Cn}z_2$.  Then
$\mu_n\nu_n$ is homotopic to $\beta^{Cn-1}$ relative to $\Fix(f)$ and
so also relative to $A$.  This proves that
$L_{\beta,\gamma}(f^n\gamma) \ge Cn-1$ and hence that $\sigma_{f,
\beta,\gamma}(\gamma) \ge C > 0$.  \endproof

We need one more definition.  
Assume that $f \in \Homeo(S)_0$ and that $S \ne S^2$.   A metric $d$ on $S$ lifts to an equivariant metric $\ti d$ on the universal cover $\ti S$.  We say that $f$ has {\it linear displacement} if either of the following conditions hold.
\begin{itemize}
\item $S \ne T^2$, $\ti f$ is the identity lift  and there exists $\ti x \in \ti S = H$ such that 
$$
\liminf_{n \to \infty} \frac{\ti d(\ti f^n(\ti x),\ti x)}{n} > 0.
$$
\item $S = T^2$ and there exist $\ti f$ and  $\ti x_1,\ti x_2 \in \ti S = \R^2$ such that
$$
\liminf_{n \to \infty} \frac{\ti d(\ti f^n(\ti x_1),\ti f^n(\ti x_2))}{n} > 0.
$$
\end{itemize}
We recognize linear displacement via the following lemma.

\begin{lemma} \label{displacement} Assume that $S \ne S^2$.  Then  $f \in \Homeo(S)_0$  has linear displacement if any of the following conditions is satisfied.
\begin{itemize}
\item $S \ne T^2$ and the identity lift  does not fix the full pre-image of $\Fix(f)$. 
\item  $S = T^2$ and there is a lift $\ti f$ that  fixes some but not every point in the full pre-image of $\Fix(f)$. 
\item  There exist  $x_1,x_2 \in S$ such that $\rho(x, \ti f) \ne \rho(x, \ti f)$ where $\ti f$ is the identity lift if $S \ne T^2$ and is arbitrary if $S = T^2$.
\end{itemize}
\end{lemma}

\proof  In the first case, there exists $\ti x \in \ti S$ and a non-trivial covering translation $T$ such that $\ti f(\ti x) = T (\ti x)$ where $\ti f$ is the identity lift.  Since $\ti f$ commutes with $T$, $\ti f^n(\ti x) = T^n (\ti x)$ for all $n$ and the lemma follows from the fact that the translation length of $T^n$ is $n$ times the translation length of $T$.  In the second case, there is a fixed point $\ti x_1$ for $\ti f$ and there is a point $\ti x_2$ such that $\ti f(\ti x_2) = T (\ti x_2)$ for some non-trivial covering translation $T$.  Thus $\ti f^n(\ti x_1) = \ti x_1$ and  $\ti f^n(\ti x_2) = T^n (\ti x_2)$ and the lemma follows as in the previous case.  In the third case the lemma follows from the definition of rotation vector. 
\endproof

The following corollary of Lemma~\ref{fromPeriodic} is the main result of this section.  

\begin{cor} \label{three possibilities} Suppose that $f \in \Diff_{\mu}(S)_0$ has infinite order and that  $\supp(\mu) \not \subset \Fix(f)$.    Suppose further that if $S = T^2$  then $\Fix(f) \ne \emptyset$ and that if $S= S^2$ then $\Fix(f)$ contains at least three points.  Then, after possibly replacing $f$ with an iterate,   at least one of the following holds:
\begin{description}
\item [(1)]There is a closed curve $\tau$ such that $\egr(f,\tau) > 0$.
\item [(2)] $f$ has linear displacement.
\item [(3)] There is a $k$-fold cover $S_k$ of $S$ with $k = 1$ or $k = 2$,  a lift $f_k : S_k \to S_k$ of $f : S \to S$ that is isotopic to the identity and there exist $\alpha,\beta$ and $\gamma$ as above such that  $\sigma_{f_k, \beta,\gamma}(\alpha) > 0$.
\end{description}
\end{cor}

\proof We consider the five cases of Lemma~\ref{fromPeriodic}. If the first item of that lemma holds then the first item of this corollary holds.  If the fifth item of the lemma holds then the second item of the corollary holds by Lemma~\ref{displacement}.  

Suppose that the third item of the lemma holds.  Let $\beta$ be the
core curve of $U$. Let $\gamma$ be an arc with interior in $U$ that
extends to an arc $\bar \gamma$ in $\bar U$ with endpoints on
different components of $\partial \bar U$ and that intersects $\beta$
exactly once.  Finally, let $\alpha$ be an arc whose interior is
contained in $U$, that has $x$ as one endpoint and that extends to an
arc $\bar \alpha$ that is disjoint from $\bar \gamma$ and that has an
endpoint in a component of $\partial \bar U$ consisting of fixed
points of $f$.  Then $\sigma_{f, \beta,\gamma}(\alpha) > 0$ and (3) is
satisfied with $k =1$.

     Suppose next that the second item of the lemma holds, i.e. that $f$ is isotopic relative to $\Fix(f)$ to a composition of Dehn twists about a finite collection $R(f)$ of simple closed curves in $S \setminus \Fix(f)$ that are neither peripheral nor parallel in $S \setminus \Fix(f)$. Denote the set of elements of $R(f)$ that are essential in $S$ by $R_e(f)$.  If $R_e(f) \ne \emptyset$  then $S \ne S^2$ and there are two possibilities.  One is that there exist $x_1, x_2 \in \Fix(f)$ and a path $\tau$ connecting them such that $\tau$ has exactly one intersection with $R_e(f)$.  In this case no lift of $f$ fixes the full pre-image of  $x_1$ and $x_2$ so (2) holds by Lemma~\ref{displacement}.  If $S = T^2$ then this case always occurs so we may now assume that $S$ has negative Euler characteristic.
 The second possibility is that there exists $x \in \Fix(f)$, a lift $\ti x \in H$ and a ray $\ti \sigma$ connecting $\ti x$ to a point in $\sinfty$ such that $\ti \sigma$ crosses exactly one element of the full pre-image $\ti R_e$ of $ R_e$.  In this case $\ti x$ is not fixed by the identity lift $\ti f$ and so again (2) holds.  

Suppose then every element of $R(f)$ bounds a disk in $S$ that contains at least two points in $\Fix(f)$. We will assume at first that $\Fix(f)$ is not contained in a single component of $S \setminus R(f)$.   Choose a path $\gamma$ connecting $x_1 \in \Fix(f)$ to $x_2 \in \Fix(f)$ such that $\gamma$ has either one or two intersections with $R(f)$ and that if there are two intersections then they are with different elements of $R(f)$. Choose a lift $\ti \gamma \subset H \cup \sinfty$ of $\gamma$ and let $\ti f : H \cup \sinfty \to H \cup \sinfty$ be the lift that fixes the initial endpoint $P \in \sinfty$ of $\ti \gamma$.  If $\beta$ is the first element of $R(f)$ intersected by $\gamma$ then there is a lift $\ti \beta$ of $\beta$ that intersects $\ti \gamma$; the endpoints of $\ti \beta$ are fixed by $\ti f$.  Let $T : H \cup \sinfty \to H \cup \sinfty$ be the covering translation corresponding to  $\ti \beta$.  Let $l$ be the order of the Dehn twist of $f$ about $\beta$ and let $Q \in \sinfty$ be the terminal endpoint of $\ti \gamma$.   If $\beta$ is the only element of $R(f)$ that intersects $\gamma$ then $\ti f^n(Q) = T^{nl}(Q)$.   If $\gamma$ intersects another element $\ti \beta'$ of $\ti R(\ti f)$, with endpoints say $R_1, R_2 \in \sinfty$, then $\ti f^n(R_i) = T^{nl}(R_i)$ and $\ti f^n(Q)$ lies in the interval bounded by $T^{nl}(R_1)$ and $T^{nl}(R_2)$ that does not contain $P$.  In either case $\ti f^n(\ti \gamma)$ intersects $T(\ti \gamma)$ and $T^{nl-1}(\ti \gamma)$.  An argument exactly like that in the proof of Lemma~\ref{spread and tracing} shows that $\sigma_{f, \beta,\gamma}(\gamma) > 0$ and (3) is satisfied with $k =1$.

To complete this case, suppose $R(f)$ is a single element $\beta$ and that $\Fix(f)$ is contained in a disk bounded by $\beta$. Then $S \ne S^2$ and we can choose  a connected two-fold cover $S_2$  of $S$ and  a lift  $f_2 :S_2 \to S_2$  of $f: S \to S$ that fixes the full pre-image of $\Fix(f)$.   The pre-images $\beta_1$ and $\beta_2$ of $\beta$ in $S_2$  bound disjoint disks  that contain at least two points in $\Fix(f_2)$.  The isotopy, relative to $\Fix(f)$, of $f$ to a  Dehn twist about $\beta$ lifts to an isotopy, relative to $\Fix(f_2)$,  of $f_2$ to a composition of Dehn twists about $\beta_1$ and $\beta_2$.  We are now reduced to the previous subcase.

It remains to consider the case that the fourth item of Lemma~\ref{fromPeriodic} holds.  After passing to a two fold cover if necessary, the hypothesis of Lemma~\ref{spread and tracing} applies so (3) is satisfied.  
\endproof

\section{Proof of Theorem~\ref{mainDistort}}  

We prove Theorem~\ref{mainDistort} by showing that $f$ can not satisfy any of the three possibilities listed in  Corollary~\ref{three possibilities}.   We begin this section by proving the lemmas that relate to (1) and (2) of that corollary .  We then state a  proposition needed for the third and prove the theorem.  The section concludes with a proof of the proposition.   

\begin{lemma} \label{no linear displacement} If $G$ is a finitely generated subgroup of $\Homeo(S)_0$ and $f \in G$ is distorted in $G$   then $f$ does not have linear displacement. 
\end{lemma}

\proof Suppose at first that $S$ has negative Euler characteristic.
In this case the identity lifts $\{\ti g: g \in G\}$ form a subgroup
$\ti G$ and $\ti f$ is a distortion element in $\ti G$.  Let $d$ be
the distance function of a Riemannian metric on $S$ and let $\ti d$ be
its lift to $H$.  For generators $g_1,\dots,g_j$ of $G$ there exists
$C > 0$ such that $\ti d(\ti g_i(\ti x),\ti x) < C$ for all $\ti x \in
H$ and all $i$.  It follows that
$$ \liminf_{n \to \infty} \frac{\ti d(\ti f^n(\ti x),\ti x)}{n} \le
\liminf_{n \to \infty} C \frac{|f^n|}{n} = 0.
$$

We now turn to the case that $S=T^2$ and $\ti S = \R^2$. 
Suppose that $\ti f : \R^2 \to \R^2$ is a lift  of $f$ and that  $\ti x_1, \ti x_2 \in \R^2$.  Choose  lifts $\ti g_i$ of the generators $g_i$ and note that there is a
constant $C> 0$ such that $\ti d(\ti g_i(\ti x), \ti x)  \le C$ for all $i$ and
all $\ti x \in \R^2$ and such that  $ \ti d(\ti x_1, \ti x_2) < C$. 
Now $f^n = w_N \dots w_2 w_1$ where $N = |f^n|$
and each $w_j$ is one of the $g_i$ or its inverse.  Let $\ti w_j = \ti
g_i$ or $\ti g_i^{-1}$ if $w_j = g_i$ or $ g_i^{-1}$ and let $\ti w =
\ti w_N \dots \ti w_2 \ti w_1$.  Then $\ti d(\ti w(\ti x),\ti x) \le CN$
for all $\ti x \in \ti \R^2$ which implies that  $\ti d(\ti w(\ti x_1),\ti w(\ti x_2)) <C(2N+1)$.  
Since $\ti w$ is a lift of $f^n,\ T \ti w = \ti f^n$ 
for some covering translation $T$.  The metric $\ti d$ is $T$-equivariant.   Thus $\ti d(\ti f^n(\ti x_1),\ti f^n(\ti x_2)) = \ti d(\ti w(\ti x_1),\ti w(\ti x_2))$ and 
\[ 
  \liminf_{n \to \infty} \frac{\ti d(\ti f^n (\ti x_1),\ti f^n(\ti x_2))}{n}
\le C (1+ 2 \liminf_{n \to \infty} \frac{|f^n|}{n} )= 0
\]
since $f$ is a distortion element.  This proves that $f$ does not have linear displacement.  
\endproof

\begin{lemma} \label{lem: distortion_per} 
If $G \subset \Homeo_{\mu}(T^2)_0$ is a finitely generated 
subgroup and $f$ is a distortion element of $G$ then 
$\Per(f) \ne \emptyset$.
\end{lemma}

\proof Since $T^2$ has no distortion elements (see 
Remark~\ref{rem: distortion}), the image of $f$ under
the homomorphism $\hat \rho_{\mu} : G \to T^2$ has finite order,
say $k$.  Let $h = f^k$. Remark~\ref{zero mean rotation} implies that
there is a lift $\ti h$ with $\rho_{\mu}(\ti h) = 0$.
Lemma~\ref{displacement} and Lemma~\ref{no linear displacement} imply
that $\rho(x,\ti h)$ is independent of $x$.  Thus $\rho(x,\ti h) = 0$
whenever it is defined and Lemma~\ref{ergodic measure} completes the proof. \endproof

\begin{lemma} \label{egr} If $G$ is a finitely generated subgroup of $\Diff(S)_0$ and $f \in G$ is distorted in $G$   then $\egr(f, \tau) = 0$ for all closed curves $\tau$.
\end{lemma}

\proof. Choose generators $g_1,\dots,g_j$ of $G$.  There exists  $C > 0$ such that $||Dg_i|| < C$ for all $i$.  Thus $l_S(g_i(\tau)) \le C l_S(\tau)$ for all $\tau$ and all $i$.  It follows that  
$$
 \liminf_{n \to \infty}\frac{\log(l_S(f^n(\tau))}{n} \le \liminf_{n \to \infty}\frac{\log(l_S(\tau)) + \log(C) |f^n|}{n} = 0.
$$
\endproof 

The following lemma is used to reduce the third item of Corollary~\ref{three possibilities}  to the case that $k =1$.

\begin{lemma} \label{no harm in lifting} Suppose that  $f \in \Diff(S)_0$, that $S_2$ is a two fold cover of $S$ and that   $f_2 \in \Diff(S_2)_0$ is a lift of $f$.   If $f$ is distorted in a finitely generated subgroup $G$ of $\Diff(S)_0$ then $f_2$ is distorted in a finitely generated subgroup $G_2$ of $\Diff(S_2)_0$.
\end{lemma}

\proof  Choose generators $g_1,\dots,g_j$ of $G$ and lifts $h_i \in \Diff(S_2)_0$  of $g_i$.  Let $t_2: S_2 \to S_2$ be a generator of the (order two) group of covering
translations of  $S_2$ and define $G_2$ to be the subgroup generated by  $\{t_2, h_1, \dots h_j\}$.    If $f^n = g_1 \dots g_m$ then
$f_2^n = t_2^i h_1 \dots h_m$, where $i = 0$ or $1$, 
so $|f_2^n|$ in these generators
is at most  $|f^n| + 1$ and $f_2$ is a distortion element of $G_2$.
\endproof

\begin{prop} \label{prop: spread}  
If $G$ is a finitely generated subgroup of $\Diff(S)_0$ and $f \in G$ is distorted in $G$   then $\sigma_{f,\beta,\gamma}(\alpha) = 0$ for all $\alpha,\beta,\gamma$.
\end{prop}

Before proving Proposition~\ref{prop: spread} we use it to prove Theorem~\ref{mainDistort}. 

\vskip .1in

\noindent{\bf Proof of Theorem~\ref{mainDistort}} \ If $S = T^2$ 
let $n$ be the smallest positive integer such that $\Fix(f)
\ne \emptyset$.
If $S=S^2$ let $n$ be the smallest positive integer such that $\Fix(f)$
contains three points. And let $n=1$ if the genus of $S$ is at least two.
Lemma~\ref{egr} and Lemma~\ref{no linear
displacement} imply that $f^n$ does not satisfy either (1) or (2) of
Corollary~\ref{three possibilities}. Lemma~\ref{no harm in lifting}
and Proposition~\ref{prop: spread} imply that $f^n$ does not satisfy (3)
of Corollary~\ref{three possibilities}. Since $f^n$ does not satisfy the
conclusion of Corollary~\ref{three possibilities}, it cannot satisfy
the hypotheses.  The only unverified hypothesis is that $\supp(\mu)
\not \subset \Fix(f^n)$, so this must fail. 
\qed

\bigskip

The proof of Proposition~\ref{prop: spread} requires some preliminary lemmas. The reason that this case is more difficult than the others is that the generators $g_j$ of $G$ need not fix the endpoints of $\gamma$ even though $f$ does. To overcome this we will show that if $f^n$ can be written as a product
of $N$ generators then it can be written as a product of $N$ diffeomorphisms 
chosen from a (perhaps infinite) collection of diffeomorphisms that  fix the endpoints of $\gamma$ and that 
possess a uniform bound on their spread.

Define $D$ to be twice the diameter of $S$.  By Remark~\ref{gamma independent} we may assume that $\gamma$ is geodesic with respect to an underlying Riemannian metric on $S$ and has  length at most $D$.

\begin{lemma} \label{lem: geodesic seg}  
Suppose that $g \in \Diff(S)$ and that 
$\eta$ and $\eta'$ are smoothly embedded geodesic arcs in
$S$ with length at most $D$.  There exists a constant $C(g)$, independent of $\eta$ and
$\eta'$ such that the absolute value of the
algebraic intersection number of any subsegment
of $g(\eta)$ with $\eta'$ is less than $C(g).$
\end{lemma}
\begin{proof}
We first show that there is an $\varepsilon >0$ such that this result
holds (in fact with $C(g) = 1$) if the lengths of $\eta$ and $\eta'$
are less than $\varepsilon >0$.  

To do this we observe that $S$ can
be covered with a finite set of charts in which geodesics are 
roughly straight, for example these charts may be obtained via the
exponential map $exp: TS_{x_i} \to S$ where we consider the domain 
to be a neighborhood of $0 \in TS_{x_i}$ and we choose the points
$x_i$ so that the charts defined cover all of $S$.  By roughly straight
we will mean that if $u$ and $v$ are unit tangent vectors to a
geodesic segment lying entirely in the image of a single chart
then the inner product $(u,v)$ of $u$ and $v$ in the co-ordinates of the 
chart is $> 1/2.$ (The inner product is well defined even though $u$ and $v$ have different base points because we are working in a chart.)

We next observe that if $\varepsilon$ is
sufficiently small and $x,y \in S$ are points within $\varepsilon$
of each other then $Dg$ is roughly linear at these points.
In particular for $\varepsilon$ sufficiently small we may assume
that if $u_x$ and  $v_y$ are unit tangent vectors at $x$ and $y$ respectively,
and they satisfy $(u_x,v_y)\ > 1/2$ in one of our charts, then
$(Dg(u_x),Dg(v_y))\ > 0$ in any of our charts containing $g(x)$ and $g(y).$

Suppose now that $\beta$ and $\beta'$ are geodesic segments of length
less than $\varepsilon$ and $g(\beta) \cap \beta' \ne \emptyset.$ We
may assume that there is a single one of our charts containing both
$g(\beta)$ and $\beta'$.  If $g(\beta)$ crosses $\beta'$ twice in
succession with the same intersection number it can do so only by
going entirely around one end of $\beta'$ (looked at in local
co-ordinates of the chart).  The mean value theorem would then say
that there are two points on $g(\beta)$ where the tangents have
opposite directions.  But this would imply the existence of tangent
vectors $u_x$ and  $v_y$ to $\beta$ with $(Dg(u_x),Dg(v_y))\  < 0.$ 
Hence the result holds with $C(g) = 1$ if the lengths of $\beta$ and
$\beta'$ are less than $\varepsilon.$

The general case now follows because the lengths of $\beta$ and $\beta'$ are  are bounded by $D$.  We can divide each of $\beta$ and $\beta'$ into 
subsegments of length at most $\varepsilon$ and apply the result
above.  It follows that the absolute value of the algebraic
intersection number of any subsegment of $g(\beta)$ with $\beta'$ is
less than $(D/\varepsilon)^2$.
\end{proof}

Let $\gamma$ be a fixed oriented geodesic arc in $S$ with length at
most $D$, let $A =\{x,y\}$ be its endpoint set and let $M$ be the
surface with boundary obtained from $S \setminus A$ by blowing up $x$
and $y$.  For each ordered pair $\{x',y'\}$ of distinct points in $S$
choose once and for all, an oriented geodesic arc $\eta = \eta(x',y')$
of length at most $D$ that connects $x'$ to $y'$ and choose $h_{\eta}
\in \Diff(S)_0$ such that $h_{\eta} (\gamma) = \eta,\ h_{\eta} (x) =
x',\ h_{\eta} (y) = y'.$ There is no obstruction to doing this since
both $\gamma$ and $\eta$ are contained in disks.  If $x = x'$ and
$y=y'$ we choose $\eta =\gamma$ and $h_{\eta} = id.$

Given  $g \in \Diff(S)$  and an ordered pair $\{x',y'\}$ of distinct points in $S$, let $\eta = \eta(x',y')$, $\eta' = \eta(g(x'),g(y'))$ and note that 
$g_{x',y'} := h_{\eta'}^{-1} \circ g \circ h_\eta$
pointwise fixes $A$. The following lemma asserts that although the pairs $\{x',y'\}$ vary  over a non-compact space, the  elements of $\{g_{x',y'}\}$ exhibit uniform behavior from the point of view of spread. 

\begin{lemma} \label{lem: gamma-len}  With notation as above, the following hold for all $g \in \Diff(S)$.   
\begin{enumerate}
\item   There exists a constant $C(g)$  
such that 
\[
L_{\beta,\gamma}(g_{x',y'}(\gamma)) \le C(g) \mbox{ for all } \beta \mbox{ and all } x',y'.   
\]
\item There exists a constant $K(g)$ such that 
\[
L_{\beta,\gamma}(g_{x',y'} (\alpha)) \le L_{\beta,\gamma}(\alpha) + K(g) \mbox{ for all } \beta, \mbox{  all } \alpha \mbox{ and all }  x',y'. 
 \]
\end{enumerate}
\end{lemma}
\begin{proof}
If $ \delta$ is a subarc of $ g_{x',y'}(\gamma)$   then $h_{\eta'}(\delta)$
is a subarc of $g(h_{\eta}( \gamma)) = g(\eta).$  Also
$h_{\eta}(\delta \cap  \gamma)$
is the intersection of this subarc of $g(\eta)$ with $\eta'.$
Assertion (1) now follows from Lemma~\ref{lem: geodesic seg} and Remark~\ref{alg int}.

We next prove assertion (2) using the notation of the $S \ne S^2$ case of  section~\ref{s: spread}; the argument works equally well in the $S = S^2$ case. There exist $a$ and   $b = a +  L_{\beta,\gamma}(\alpha)$ such that $\ti \alpha$ is in the region bounded by 
 $\ti \gamma_a$ and  $\ti \gamma_b$.  Since  $g_{x',y'}$ fixes $x$ and $y$ it extends to a map (also called)  $g_{x',y'}: M \to M$.  Choose a lift $\ti g_{x',y'}:  H \cup \sinfty \to H \cup \sinfty$ such that $\ti g_{x',y'} (\ti \gamma_a)$ intersects the region bounded by $\ti \gamma_a$ and $\ti \gamma_{a+1}$.  Let $a' =  a-C(g) -1$ and $b' = b + C(g) + 1$.  By (1), the region bounded by $\ti \gamma_{a'}$ and $\ti \gamma_{b'}$ contains    $\ti g_{x',y'}(\ti \gamma_a)$ and $\ti g_{x',y'}(\ti \gamma_b)$ and hence also $\ti g_{x',y'}(\ti \alpha)$.  Thus   $L_{\beta,\gamma}(g_{x',y'} \alpha) \le b' - a' \le  L_{\beta,\gamma}(\alpha) + 2C(g) + 2$ and we may set $K(g) =  2C(g) + 2$. 

\end{proof}

The elements of $\{g_{x',y'}\}$ play a crucial role in the next lemma.

\begin{lemma} \label{T-len growth}  
Suppose that $g_i \in \Diff(S)_0,\  1 \le i\le k,$ that $f$  is in the group they generate and that   
  $|f^n|$ is the word length of $f^n$ in the generators $\{g_i\}$. Then there is a constant $C >0$ such that
$$
L_{\beta,\gamma}(f^n (\alpha)) \le L_{\beta,\gamma}(\alpha) + C |f^n| 
$$ 
for all $\alpha,\beta,\gamma$ and all $n >0.$
\end{lemma}
\begin{proof} We may assume that $\gamma$ is geodesic with length at most $D$.
For each $g_i$, let $K(g_i)$ be the constant given by Lemma~\ref{lem: gamma-len}
and let $C = max\{K(g_i)\}.$  Suppose $N = |f^n|$ and 
$f^n = w_Nw_{N-1} \dots w_2w_1$ where each $w_j$ is one of the generators
$\{g_i\}$ or its inverse.  Let $x$ and $y$ be the endpoints of $\gamma$, let  
  $\eta_j$ be the chosen  geodesic joining 
$w_j w_{j-1}\dots w_1(x)$ and $w_j w_{j-1}\dots w_1(y)$ and let
$h_j = h_{\eta_j}$ as above  be a diffeomorphism carrying $\gamma$
to $\eta_j$. Denote $h_j^{-1}w_jh_{j-1}$ by $\hat w_j$.

We note that $h_0 =  id$ and since $f^n$ fixes $x$ and $y$ we
also have $h_N =  id$.  Hence
\begin{align*}
 f^n &= w_N w_{N-1} \dots w_2 w_1\\
&= (h_N^{-1} w_N h_{N-1})(h_{N-1}^{-1} w_{N-1} h_{N-2})
\dots (h_2^{-1} w_2 h_1) (h_1^{-1} w_1 h_0),
\end{align*}

so 
\[
 f^n = \hat w_N \hat w_{N-1} \dots \hat w_2 \hat w_1.
\]

Now $N$ applications of Lemma~\ref{lem: gamma-len}
 yield
\begin{align*}
L_{\beta,\gamma}( f^n(\alpha))
&= L_{\beta,\gamma} (\hat w_N \dots \hat w_2 \hat w_1(\alpha))\\
&\le L_{\beta,\gamma} (\hat w_{N-1} \dots \hat w_2 \hat w_1(\alpha))  + C\\
&\dots\\
&\le L_{\beta,\gamma} (\alpha) +C N.
\end{align*}
\end{proof}

\noindent{\bf Proof of Proposition~\ref{prop: spread}}
Since $f$ is distorted in $G$
\[
\liminf_{n \to \infty} \frac{|f^n|}{n} = 0.
\]
According to the definition of spread and 
Lemma~\ref{T-len growth}
we then have
\[
\sigma_{f,\beta,\gamma}(\alpha) 
=  \liminf_{n \to \infty} \frac{L_{\beta,\gamma}( f^n(\alpha))
}{n}
\le \liminf_{n \to \infty} \frac{L_{\beta,\gamma}(\alpha) + C |f^n|}{n}
= 0.
\]\qed

\section{Applications of the main theorem}

We need a few  preliminary lemmas.

\begin{lemma} \label{lem: dist elt} 
Let $\G$ be a finitely generated almost simple group
and let $\H$ be a normal subgroup.  If 
$u \in \H$ is a distortion element in $\G$ then it is also a
distortion element in $\H$.
\end{lemma}
\begin{proof}
Since $\G$ is almost simple and $\H$ contains an infinite order
element, $\H$ has finite index.  Since $H$ is a finite index subgroup
of a finitely generated group, it is finitely generated.  Let
$\{h_1,\dots,h_m\}$ be a set of generators for $\H$ and let
$1=a_0,a_1,\dots, a_k$ be coset representatives of $\H$.  Define word
length in $\H$ with respect to $\{h_1,\dots,h_m\}$ and word length in
$G$ with respect to the generating set
$\{h_1,\dots,h_m,a_1,\dots,a_k\}$.

For each $g \in \G$  there exists $0 \le i \le k$ and $h \in \H$ such that $g = a_ih$.  Choose $C$ such that $|h|_{\H} \le C$ for the finitely many  $g$ satisfying $|g|_{\G}  \le 2$.   We claim that $|h|_{\H} \le C \cdot |g|_{\G}$ for all $g$. The proof is by induction on $|g|_{\G}$ with the $|g|_{\G}  \le 2$ being by construction.   Given $ g \in \G$ with $l := |g|_{\G}  \ge 3$,  write $g = u g_1$ where $|u|_{\G}  = 1$ and   $|g_1|_{\G}  = l-1$.  The inductive hypothesis implies that $g_1 = a_ih_1$ and $ua_i = a_j h_2$ where $|h_1|_{\H} \le C (l-1)$ and $|h_2|_{\H} \le C$.  Setting $h = h_2h_1$ we have $|h|_{\H} \le C\cdot l$ which completes the inductive step.    

     It follows that $|h|_{\H} \le C \cdot |h|_{\G}$ for all  $h \in \H$ which implies the lemma.
\end{proof}

\begin{lemma} \label{lem: no R hom} 
Let $\G$ be a finitely generated almost simple group which contains
a distortion element and let $\H$ be a normal subgroup.  Then
any homomorphism $\psi: \H \to \R$ is trivial.
\end{lemma}

\begin{proof}
Since $\G$ is almost simple, $\H$ is either finite or has finite index.
Clearly the result is true if $\H$ is finite, so we assume it has
finite index. If $u$ is a distortion element in $\G$ then $v:=u^k \in
\H$ for some $k > 0$.  The normalizer $\D$ of $v$ in $\G$ is infinite
and hence has finite index in $\G$; it is obviously contained in $\H$.  Thus
$\D$ has finite index in $\H$.  Since $\R$ contains neither torsion
nor distortion elements, $v$, and hence $\D$ is in the kernel of
$\psi$ for every homomorphism $\psi: \H \to \R$.  Since $\D$ has
finite index in $\H$ we conclude that $\psi(\H)$ is finite and hence
trivial.
\end{proof}

\begin{lemma} \label{lem: dist_heis} The generator $f$ of the center of the three dimensional Heisenberg group $\H$ is a distortion element.  
If $\psi: \H \to GL(2,\R)$ is a homomorphism then $\psi(f) = \pm I.$
\end{lemma}

\begin{proof} Suppose $\H$ is the three dimensional Heisenberg group.
There are generators $g$ and $h$ of $\H$ such that $f = [g,h]$
has infinite order and generates the center of $\H$.  It follows
that $f^{n^2} = [g^n,g^n]$ and hence that $f$ is a distortion
element of $H$. 

Suppose $\phi: \H \to GL(2,\R)$ is a homomorphism.
Let $F = \phi(f), \ G = \phi(g),\ H = \phi(h)$, all of
which we consider as $2 \times 2$ matrices.

If $F$ has distinct eigenvalues then $F, G$ and
$H$ are all simultaneously diagonalizable over $\mathbb C$.
Hence all commute and it follows that $F = [G,H] = I$.

Suppose $F$ has a single eigenvalue.  Since $F$ is a commutator
its determinant is $1$ and hence its eigenvalue must be $1$ or $-1$.
If $F$ is diagonal we are done so we may assume  $F$ is
upper triangular in a suitable basis with a non-zero off-diagonal
entry.  An easy computation shows
the fact that $F$ commutes with $G$ and $H$ implies that they
also are upper triangular in this basis and have equal diagonal
entries.  From this one sees again that $G$ and $H$ commute, so again $F=I$.
\end{proof}

\begin{lemma}\label{lem: heis_per}
Suppose the subgroup $H \subset \Diff_{\mu}(S^2)_0$ is finitely generated and that $f$ is a central element of $H$ that is distorted in $H$.  Suppose further that the support of $\mu$ contains at least three points.  Then $\Per(f)$ contains at least three points.  
\end{lemma}

\proof  Since $f$ is orientation preserving, it has at least one fixed
point.  If there were only one fixed point then its complement would
be an invariant open disk $U$ with $\mu(U) > 0$.  Poincar\'e recurrence
implies that $U$ contains recurrent points for $f$ and the Brouwer
translation theorem then implies that $U$ contains fixed points for
$f$.  Thus $\Fix(f)$ contains at least two points.

Suppose that $\Fix(f)$ has exactly two points.  
Since $f$ is central, every element of $H$ preserves $\Fix(f)$.
If the elements of $\Fix(f)$ are atoms for the measure $\mu$ we change $\mu$
so $\mu(\Fix(f)) = 0.$  Let $A$ be the closed
annulus obtained by blowing up these points.  There is an induced
measure, still called $\mu$, such that $\mu(\partial A) = 0$.  
After passing to an index two subgroup of $H$, we may assume that every
element of $H$ fixes the elements of $\Fix(f)$ and so $H$ can be
thought of as a group of homeomorphisms of $A$ that preserve
orientation and the components of $\partial A$.  Consider the
homomorphism $\hat \rho_{\mu} : H \to S^1$ of section~\ref{rotation
vectors}.  Since $S^1$ has no distortion elements, some iterate of
$f$, say $f_1 = f^m$, is in the kernel of this
homomorphism. Remark~\ref{zero mean rotation} implies that the mean
rotation vector $\rho_{\mu}(\ti f_1) = 0$ for some lift $\ti f_1$.

Letting each element of $H$ act on the double of $A$ defines an
embedding of $H$ into $\Homeo_\nu(T^2)$, where the measure $\nu$ on $T^2$   is the  double of the measure $\mu$ on $A$.   
Lemma~\ref{displacement} and Lemma~\ref{no linear displacement}, applied to the double of $f_1$, imply  
that $\rho(x,\ti f_1)$ is independent of $x$.  Thus $\rho(x,\ti f_1) =
0$ whenever it is defined and $\rho_{\nu}(\ti f_1) = 0$ for all
$f_1$-invariant measures $\nu$.  Choosing $\nu$ to be an ergodic
invariant measure the result follows from Proposition~(\ref{prop:
ergodic_annulus}).  \endproof

We say that $x \in S$ is a {\em global fixed point} for $H$ if $x
\in \Fix(h)$ for each $h \in H$.

\begin{lemma}\label{lem: global_fp}
Suppose $G$ is a finitely generated subgroup of $\Diff(S)$ which is
almost simple and contains a distortion element $f$.  If $H
\subset G$ is a normal subgroup that contains $f$ then the set  
$\Lambda$  of global fixed points for $H$  is finite.
\end{lemma}
\begin{proof}
The subgroup $H$ is normal so by 
Lemma~(\ref{lem: no R hom}) the only homomorphism 
from $H$ to $\R$ is the trivial one.  We will use this fact repeatedly.

We will assume $\Lambda$ is infinite and show this leads to a
contradiction.  
Let $x \in \Lambda$
be an accumulation point and let $v \in TS_x$ be a
limit of a subsequence from the sequence of vectors $\{v_n\}$ 
where $v_n = (x - x_n)/\|x - x_n\|$ (in local co-ordinates) 
and $x_n \in \Lambda$
satisfies $\displaystyle{ \lim_{n\to \infty} x_n = x.}$ Then
$Dh_x(v) = v$ for all $h \in H$ and $v$ is a global fixed point
of the action of $H$ on the circle $S^1$ of unit
vectors in $TS_x$ where $h$ acts on $S^1$ as the projectivization of 
$Dh_x$ on $TS_x$.
Considering this action of $H$ on $S^1$ and taking the
logarithm of the derivative at the global fixed point $v$ gives a 
homomorphism from $H$ to $\R$,
which must be trivial.

Thus the value of these derivatives at
$v \in S^1$ is $1$ for every element of $H$, so we can apply the Thurston
stability theorem (Theorem 3 of \cite{Th}) to this action.
This theorem asserts that either the
action of $H$ on $S^1$ is trivial or there is a non-trivial
homomorphism from $H$ to $\R.$  Since the latter is impossible
we conclude that the action of $H$ on $S^1$ is trivial. In fact
since the group is acting on $S^1$ by fractional linear transformations
it is not difficult to give a direct argument that $H$ acts trivially.

It follows that for every $h \in H$ we have $Dh_x = I$,
since $Dh_x$ must be a homothety to projectivize to the identity,
but also  $Dh_x$ fixes the vector $v$.

We can now again apply the Thurston stability theorem, this time to $H$ acting on $S$ with
the global fixed point $x$.  The theorem asserts that
either $H$ acts trivially or there is a non-trivial
homomorphism from $H$ to $\R$.  Since neither of these are
the case we have arrived at a contradiction.
\end{proof}

\noindent{\bf Theorem~\ref{thm: finite_image}} {\em 
Suppose that $S$ is a closed oriented surface of positive genus
equipped with a Borel probability measure $\mu$ 
and that $\G$ is a finitely generated group which is almost simple and
possesses a distortion element $u$.  Suppose further that either $\mu$ has infinite support or that $\G$ is a Kazhdan group.  Then any homomorphism $\phi: \G
\to \Diff_{\mu}(S)$ has finite image.  The result is valid for $S =
S^2$ with the additional assumption that $\phi(u)$ has at least three
periodic points. }  

\begin{proof}
Since $\G$ is almost simple, it suffices to show that the kernel
of $\phi$ is infinite.  We assume that the
kernel of $\phi$ is finite and argue to a contradiction.  

Since $f = \phi(u)$ has infinite
order in $G = \phi(\G)$, it is a distortion element in $G$.   By \cite{flm}, there are no distortion elements in $\MCG(S)$.  After replacing $f$ with an iterate we may assume that $f \in \Diff_\mu(S)_0.$  Lemma~(\ref{lem: dist elt}) implies that $f$ is a distortion element in $\Diff_\mu(S)_0.$  

In the case that $S = T^2$ we note that Lemma~\ref{lem: distortion_per} 
that $\Per(f) \ne \emptyset.$
After replacing $f$ with a further iterate if necessary, we may assume by  
Theorem~(\ref{mainDistort})  that $\supp(\mu) \subset \Fix(f).$  Since $\supp(\mu)$ is invariant under $G$,  the stabilizer  
$H = \{ g \in G\ |\ g(x) = x \text{ for all } x \in
\supp(\mu)\}$ of $\supp(\mu)$ is a normal subgroup of $G$.  Lemma~\ref{lem: global_fp}
implies that $\supp(\mu)$ is finite. It follows from our hypotheses that $\G$ is a Kazhdan group.   A result of R. Zimmer
(Theorem 3.14 of \cite{Z3})   implies that every $\G$ orbit is infinite.  But this contradicts the fact that each point in $\supp(\mu)$ has finite $H$-orbit and hence finite $G$-orbit  since $H$  has finite index in $G$.

\end{proof}

\noindent{\bf Theorem~\ref{thm: Heisenberg}} 
{\em 
 Suppose that $S$ is a closed oriented surface
equipped with a Borel probability measure $\mu$, and that $\G$ is a
finitely generated almost simple group which has a subgroup $\H$
isomorphic to the three-dimensional Heisenberg group.
Then any homomorphism $\phi: \G \to \Diff_{\mu}(S)$ has finite image.}

\begin{proof}
Since $\G$ contains a subgroup isomorphic to the Heisenberg group
we know by Lemma~(\ref{lem: dist_heis}) that it contains a distortion
element.   If $\supp(\mu)$ is infinite,   then the theorem follows from 
Theorem~\ref{thm: finite_image} with Lemma~(\ref{lem: dist_heis}) and Lemma~(\ref{lem: heis_per}) supplying
the extra hypothesis needed when $S = S^2.$

Suppose  now that $\supp(\mu)$ is finite.   Then the stabilizer
$G_1 = \{ g \in G\ |\ g(x) = x \text{ for all } x \in
\supp(\mu)\}$ of $\supp(\mu)$ is an infinite, and hence finite index, normal subgroup of $G$.  Choose $x \in \supp(\mu)$. After replacing $\H$ with a finite index subgroup of itself, we may assume that $\H \subset \G_1$.  

Choose $x \in \supp(\mu)$.  The group $\G_1$ acts via the derivative on $TS_x$  and hence there is a homomorphism 
$\psi_x: \G_1 \to GL(2,\R)$ defined by $\psi_x(g) = D\phi(g)_x.$
If $f$ is the generator of the center of $\H$ then 
Lemma~(\ref{lem: dist_heis}) implies $\psi_x(f^2) = I$. We conclude that the kernel $\G_2$ of $\psi_x$ is infinite and hence has finite index in $\G$.

The Thurston stability theorem (Theorem 3 of \cite{Th}) applied
  to $\G_2$ acting with global fixed point $x$ implies that
  either $\phi(\G_2)$ is trivial
or there is a non-trivial homomorphism from $\G_2$ to $\R.$
But Lemma~(\ref{lem: no R hom}) asserts there is only the
trivial homomorphism from $\G_2$ to $\R$ so $\phi(\G_2)$ is trivial and $\phi(\G)$ is finite.  
\end{proof}

\noindent{\bf Proof of Corollary~\ref{cor: lattice}} \  By
\cite{lmr}, $\G$ has a distortion element $u$.  The Margulis normal
subgroups theorem (see \cite{Mar}) implies that $\G$ is almost simple and 
results of \cite{K} imply that $\G$ is a Kazhdan group.
 Thus Theorem~\ref{thm: finite_image}
implies the desired result.  \qed

\bigskip

The following easy and well known lemma is needed in the proofs of
Corollary~\ref{cor: commutator} and Corollary~\ref{cor: nilpotent}.

\begin{lemma} \label{finite order} 
Suppose that $f \in \Diff(S)_0$ has finite order.  If $S = T^2$ assume
that $\Fix(f) \ne \emptyset$ and if $S = S^2$ assume that $\Fix(f)$
contains at least three points.  Then $f$ is the identity.
\end{lemma}

\proof By averaging a Riemannian metric we may assume $f$ is an
isometry. If $f$ is not the identity then all fixed points of $f$ are
isolated and have index $+1$.  By hypothesis the number of fixed
points of $f$ is greater than the Euler characteristic of $S$.  This
contradicts the Lefschetz fixed point theorem unless $f = id.$
\endproof

\noindent{\bf Proof of Corollary~\ref{cor: commutator}} 
Suppose first that the genus of $S$ is greater than one.
Then Theorem~(\ref{mainDistort})
implies that $S = \supp(\mu) \subset \Fix(g)$ for every 
$g \in [\phi(\G), \phi(\G)].$
Thus $[\G, \G]$ is in the kernel of $\phi$ and $\phi(\G)$ is abelian.  

In case $S = T^2$ Lemma~\ref{lem: distortion_per} 
says that $\Per(f) \ne \emptyset.$
We know that $[\phi(\G), \phi(\G)]$ is generated
by the $\phi$ image of distortion elements, i.e. according to
Theorem~(\ref{mainDistort}) by elements of finite order.
Let $g$ be a finite order element of a generating set for 
$[\phi(\G), \phi(\G)]$.  Since $[\phi(\G), \phi(\G)]$ is also 
generated by commutators
the rotation vector (see section~\ref{rotation vectors}) 
$\hat \rho_\mu(g) = 0.$

Hence there is a lift $\ti g$ to $\R^2$ with mean rotation vector
$\rho_\mu(\ti g) = 0.$ This implies $\rho_\mu(\ti g^k) = 0$ for all
$k$ and hence if $k$ is the order of $g,$ then $\ti g^k$ is a covering
translation with zero mean rotation and hence it is the identity.
Thus $\ti g$ has a fixed point and $g$ must also. Lemma~\ref{finite order} implies that $g$ is the identity.

Since this is true for $g$ in a generating set of $[\phi(\G),
\phi(\G)]$ we again conclude that $[\G, \G]$ is in the kernel of
$\phi$ and that $\phi(\G)$ is abelian.   \qed

\bigskip

\noindent{\bf Proof of Corollary~\ref{cor: nilpotent}}
Suppose 
\[
\N = \N_1 \supset \N_2 \supset \dots \supset \N_m \supset \{1\}
\]
is the lower central series of the nilpotent group $\N.$  We will
show that (perhaps after passing to an index two subgroup if 
$S = S^2$) the assumption that $m > 1$ leads to a contradiction.

Assume that $m > 1$.  Then $\N_{m}$ is contained in the center of $\N$
and is generated by elements $f $ that are commutators, $f =[g,h]$, in
$\N$. We claim that $f$ has finite order.  If not, then $g$ and $h$
generate a group $\H \subset \N$ isomorphic to the Heisenberg group
and $f$ is a distortion element by Lemma~(\ref{lem: dist_heis}).  In
the case that $S = S^2$, Lemma~(\ref{lem: heis_per}) implies that $f$
has at least three periodic points and in the case $S = T^2$
Lemma~\ref{lem: distortion_per} implies that $\Per(f) \ne \emptyset.$
Thus for any genus, we have a
contradiction to Theorem~(\ref{mainDistort}) and the fact that
$\supp(\mu) = S$. This verifies the claim.

   If  the genus of $S$ is greater than $1$, then $f=id$ by Lemma~\ref{finite order}. Since this is true for  a generating set for $\N_m$, we conclude that $\N_m$ is trivial
contradicting the assumption $m>1.$

In the case $S = T^2$, since $f$ is a commutator, there is a lift $\ti f$
to $\R^2$ with mean rotation vector $\rho_\mu(\ti f) = 0.$  If $k$ is the order of $f$, then $\ti f^k$ is a covering translation with  $\rho_\mu(\ti f^k) = 0$ and so $\ti f^k = id$.   It follows that $\ti f$, and hence $f$, has a fixed point.  Thus Lemma~\ref{finite order} applies and the proof concludes as in the previous case. 

It remains to consider the case $S = S^2$.  As a special case, assume that there are a pair of global fixed points for $\N$ and let $A$ be the closed annulus obtained by blowing up these two points. Then $\N$ can be thought of as a group of homeomorphisms of $A$
that preserve orientation and the components of $\partial A$.
Consider the homomorphism $\rho_{\mu} : G \to S^1$ of
section~\ref{rotation vectors}.  Since $f$ is a commutator $\hat
\rho_{\mu}(f) = 0.$  There is a lift $\ti f$ to the universal
covering space of $A$ such that $\rho_{\mu}(\ti f) = 0.$  The proof in this  special case now concludes as in the $S = T^2$ case.

We now return to the general $S = S^2$ case.  Poincar\'e recurrence and the Brouwer plane
translation theorem imply that $\Fix(f)$ contains at least two points.  If $\Fix(f)$ has at least three points then Lemma~\ref{finite order} applies and the proof concludes as before.  Suppose then that $\Fix(f)$ has exactly two points.  Since $f$ is central, the subgroup $G$ of $\N$ consisting of elements that fix both elements of $\Fix(f)$ has index at most two.  By the special case, $G$ is abelian and we are done.  
\qed


\begin{thebibliography}{100}

\bibitem{a}
{\bf E.~ Alibegovi\'c,}
\newblock{ \em Translation lengths in $Out(F_n)$,} 
\newblock Geom. Dedicata, {\bf 92}
(2002), 87--93


\bibitem{brn-kis}
{\bf J. ~Kister and M.~Brown,}
\newblock {\em Invariance of complementary domains of a fixed point set,}
\newblock Proc. Amer. Math. Soc. {\bf 91} no.3, (1984)  503--504

\bibitem{BM}
{\bf M.~Burger and N.~Monod}
\newblock {\em Bounded cohomology of lattices in higher rank Lie groups}
\newblock J. Eur. Math. Soc. {\bf 1} (1999) 199--­235

\bibitem{FM}
{\bf B.~Farb and H.~Masur,}
\newblock {\em Superrigidity and mapping class groups,}
\newblock {Topology} {\bf 37} {no. 6} (1998) 1169--1176

\bibitem{flm}
{\bf B. ~Farb, A. ~Lubotzky and Y.~Minsky}
\newblock {\em Rank 1 phenomena for mapping class groups}
\newblock Duke Math. J. {\bf 106} no. 3 (2001)  581--597


\bibitem{FS}
{\bf B.~Farb and P.~Shalen,}
\newblock {\em  Real analytic actions of lattices,}
\newblock {Invent. Math.} {\bf 135}  (1998),  271--296. 

\bibitem{Fr}
{\bf J.~Franks,}
\newblock {\em Recurrence and Fixed Points of Surface Homeomorphisms,}
Ergodic Theory and Dynamical Systems {\bf 8*} (1988) 99--107

\bibitem{Fr5}
{\bf J.~Franks,}
\newblock {\em  Rotation vectors and fixed points of area preserving surface diffeomorphisms,}
\newblock {Trans. Amer. Math. Soc.} {\bf 348} (1996) 2637--2662


\bibitem{Fr6}
{\bf J.~Franks,}
\newblock {\em  Rotation numbers for Area Preserving Homeomorphisms of the Open Annulus,}
\newblock {Proceedings of the International Conference Dynamical 
Systems and Related Topics,} K Shiraiwa, ed. 
\newblock World Scientific (1991) 123--128

\bibitem{fh:periodic}
{\bf J.~Franks and M.~Handel,}
\newblock {\em  Periodic points of Hamiltonian surface diffeomorphisms,}
\newblock  {Geom. Topol.} {\bf 7} (2003) 713--756

\bibitem{fh:group}
{\bf J.~Franks and M.~Handel,}
\newblock {\em  Area preserving group actions on surfaces,}
\newblock Geom. Topol. {\bf 7} (2003) 757--771


\bibitem{gs}
{\bf S. Gersten and H. Short,}
\newblock {\em Rational subgroups of biautomatic groups,}
\newblock {Annals of Math.} {\bf 134} (1991) 125--158

\bibitem{G}
{\bf \'E.~Ghys,}
\newblock{\em Sur les groupes engendr\'e par des 
diff\'eomorphismes proches de l'identit\'e,}
\newblock {Bol. Soc. Brasil. Mat. (N.S.)}  {\bf 24} (1993) 137­178

\bibitem{han:fpt}
{\bf M.~Handel,}
\newblock {\em    A fixed point theorem for planar homeomorphisms,}
\newblock {Topology} {\bf 38} (1999) 235--264

\bibitem{K}
{\bf D.~Kazhdan,}
\newblock {\em  Connection of the dual space of a group with the structure
of its closed subgroups,}
\newblock {Funct. Anal. Appl.} {\bf 1} (1967) 63-65


\bibitem{lmr}
{\bf A. Lubotzky, S. Mozes and M.S. Raghunathan,}
{\em The word and Riemannian metric on lattices in semisimple Lie groups,}
IHES PUbl. Math. 91 (2000), 5--53


\bibitem{Mar}
{\bf G.A. Margulis,}
\newblock {\em  Discrete subgroups of semisimple Lie groups,}
\newblock {Ergebnisse der Mathematik und ihrer Grenzgebiete}, {\bf 17.}
\newblock Springer-Verlag, Berlin (1991)

\bibitem{P}
{\bf L.~Polterovich,}
\newblock {\em  Growth of maps, distortion in groups and symplectic geometry,}
\newblock  Invent. Math.  {\bf 150}  (2002), 655--686

\bibitem{Th}
{\bf W. Thurston}
\newblock{A generalization of the Reeb stability theorem}
\newblock  Topology {\bf 13} (1974) 347­352

\bibitem{Z1}
{\bf R.~Zimmer,}
\newblock {\em  Ergodic Theory and Semisimple Groups,}
\newblock Monographs in Math., {\bf 81} Birhauser, (1984)

\bibitem{Z2}
{\bf R.~Zimmer,}
\newblock {\em  Actions of semisimple groups and discrete subgroups,}
\newblock Proc. Internat. Congr. Math. (Berkeley 1986), Vol 2
\newblock A. W. Gleason, ed., Amer. Math. Soc., Providence, RI, (1987) 1247--1258

\bibitem{Z3}
{\bf R.~Zimmer,}
\newblock {\em  Lattices in semisimple groups 
and invariant geometric structures on compact manifolds,}
\newblock in Discrete Groups in Geometry and Analysis,
\newblock Progress in Mathematics {\bf 67} Birhauser, (1987), 152--210

\end{thebibliography}
\end{document}